\documentstyle[amssymb,12pt]{amsart}

\newtheorem{thm}{Theorem}[section] 
\newcommand{\bthm}{\begin{thm}}
\newcommand{\ethm}{\end{thm}}
\newtheorem{nnthm}{Theorem} 
\newcommand{\bnnthm}{\begin{nnthm}}
\newcommand{\ennthm}{\end{nnthm}}
\newtheorem{cor}[thm]{Corollary}
\newcommand{\bcor}{\begin{cor}}
\newcommand{\ecor}{\end{cor}}
\newtheorem{lemma}[thm]{Lemma}
\newcommand{\blemma}{\begin{lemma}}
\newcommand{\elemma}{\end{lemma}}
\newtheorem{prop}[thm]{Proposition}
\newcommand{\bprop}{\begin{prop}}
\newcommand{\eprop}{\end{prop}}
\newtheorem{factnum}[thm]{Fact}
\newcommand{\bfactnum}{\begin{factnum}}
\newcommand{\efactnum}{\end{factnum}}
\theoremstyle{definition}
\newtheorem{defin}[thm]{Definition}
\newcommand{\bdef}{\begin{defin} \em }
\newcommand{\edefin}{ \end{defin}}
\theoremstyle{remark}
\newtheorem{fact}{Fact}        
\newcommand{\bfact}{\begin{fact}}
\newcommand{\efact}{\end{fact}}
\newtheorem{claim}{Claim}   
\newcommand{\bclaim}{\begin{claim}}
\newcommand{\eclaim}{\end{claim}}
\newtheorem{question}{Question}   
\newcommand{\bquestion}{\begin{question}}
\newcommand{\equestion}{\end{question}}
\newtheorem{remark}{Remark}  
\newcommand{\bremark}{\begin{remark}}
\newcommand{\eremark}{\end{remark}}
\newtheorem{remarkss}{Remarks}  
\newcommand{\bremarkss}{\begin{remarkss}}
\newcommand{\eremarkss}{\end{remarkss}}
\newcommand{\bproof}{\begin{pf}}
\newcommand{\eproof}{\end{pf}}
\newcommand{\bproofclaim}{\begin{pf*}{Proof of Claim}}
\newcommand{\eproofclaim}{\renewcommand{\qed}{\ $\square$~Claim.}\end{pf*} }

\newcommand{\dom}{\mathrm{dom}}

\renewcommand{\a}{\;\&\;}
\newcommand{\sub}{\subseteq}
\newcommand{\elem}{\prec}
\newcommand{\concat}{ \raisebox{.7ex}{$\smallfrown$} }
\newcommand{\contains}{\supseteq}

\newcommand{\map}[2]{\!:#1 \rightarrow #2}
\newcommand{\restricted}{\restriction}
\newcommand{\applied}{''}

\newcommand{\compose}{\circ}

\newcommand{\seq}[1]{{\langle #1 \rangle} }
\newcommand{\quot}[1]{{\m{``} #1 \m{''}} }
\newcommand{\m}[1]{\mbox{#1}}

\newcommand{\bp}{{\Bbb P}}

\newcommand\Qd{\mbox{{\bf Q}$_{<\delta}$}}
\newcommand\Pd{\mbox{{\bf P}$\!_{<\delta}$}}

\newcommand{\ga}{{\alpha}}     
\newcommand{\gb}{{\beta}}      
\newcommand{\gd}{{\delta}}

\newcommand{\gee}{{\eta}}

\newcommand{\gk}{{\kappa}}  
\newcommand{\gl}{{\lambda}}    
\newcommand{\gu}{{\mu}}        
\newcommand{\gn}{{\nu}}

\newcommand{\gp}{{\pi}}        
       
\newcommand{\gs}{{\sigma}}     
\newcommand{\gt}{{\tau}}

\newcommand{\gw}{{\omega}}

\newcommand\cA{\cal A}
\newcommand\cB{\cal B}
\newcommand\cC{\cal C}

\newcommand\cF{\cal F}

\newcommand\cP{\cal P}

\newcommand\cS{\cal S}

\newcommand{\set}[2]{\{\; #1 \; \mbox{$\vert$}  \; #2 \; \} } 
\newcommand{\card}[1]{{\vert #1 \vert} }

\newcommand{\forces}{{\;\Vdash} }

\begin{document}

\title{Precipitous Towers of Normal Filters}

\author{Douglas Burke}

\address{Department of Mathematics \\
        University of North Texas  \\
        Denton, TX 76203}

\email{dburke@@unt.edu}

\maketitle

\begin{abstract}
We prove that every tower of normal filters of height $\gd$ ($\gd$ 
supercompact) is precipitous assuming that each normal filter in the
tower is the club filter restricted to a stationary set.  We give
an example to show that this assumption is necessary.
We also prove that every normal filter can be generically extended to a
well-founded $V$-ultrafilter
(assuming large cardinals).
\end{abstract}

In this paper we investigate towers of normal filters.  These towers
were first used by Woodin in \cite{proc}.
Woodin proved that if $\gd$ is a Woodin cardinal and
$\bp$ is the full stationary tower up to $\gd$
($\Pd$) or the countable version ($\Qd$) then the generic
ultrapower is closed under $<\gd$ sequences (so
the generic ultrapower is well-founded). 
 We show that if $\bp$
is a tower of height $\gd$, $\gd$ supercompact, and the filters generating
$\bp$ are the club filter restricted to a
stationary set, then $\bp$ is precipitous.  We 
give an example (assuming large cardinals)
 of a non-precipitous tower.
 We also  show that
every normal filter can be extended to a $V$-ultrafilter
with well-founded ultrapower in {\em some} generic extension of $V$ 
(assuming large cardinals).  Similarly for any tower of inaccessible
height.  This is accomplished by showing that there
is a stationary set that projects to
the filter or the tower and then forcing
with $\Pd$ below this stationary set.

An important idea in our proof of precipitousness
 (Theorem \ref{thm.semiproper}) has the
following form in Woodin's proof.  If $\cA_{i} \sub \Pd$ are 
maximal antichains ($i \in \gw$ and $\gd$ Woodin) then there
is a $\gk < \gd$ such that
$$ \left.  \begin{array}{cl}
   \{a \elem V_{\gk +1} \mid  & \card{a} < \gk \a (\forall i \in \gw) \:
			    \exists b \elem V_{\gk + 1} \m{ such that} \\
	         &       1)\;\; a \sub b, b \m{ end extends } a \cap V_{\gk} \\
		 &	    2)\;\; \exists x \in \cA_{i} \cap V_{\gk} \cap
			       b \: (b \cap \cup x \in x)
   \end{array}   \right\}
$$
contains a club (relative to $\card{a} < \gk$).

Before this, a similar idea was used in \cite{fms}.  For example, let
 $\cA \sub NS_{\gw_{1}}$
be a maximal antichain.  If the sealing off forcing for $\cA$ is semiproper
(this holds if we collapse a supercompact cardinal to $\gw_{2}$) then
$$
    \set{a \elem H(\gl)}{\card{a} < \gw_{1} \a (\exists b) \: b \cap \gw_{1}
			 = a \cap \gw_{1} \a \exists A \in \cA \cap b \:
			 (b \cap \gw_{1} \in A) }
$$
contains a club in $\cP_{\gw_{1}}(H(\gl))$ ($\gl >> \gw_{1}$).  This
can be  used to show that $NS_{\gw_{1}}$ is precipitous.

The basic facts about forcing that we use can be found in \cite{jech.book}
or \cite{kunen.book}.  Throughout this paper generic mean set generic.

The author is indebted to many people for helpful remarks and suggestions;
including M. Foreman, S. Jackson, T. Martin, and H. Woodin.  Some of the
results in this paper appeared in the author's thesis, supervised
by T. Martin.

   \section{Normal Filters}

In this section we recall  the basic definitions and
facts about   normal filters on $\cP(X)$. 
The proofs of these facts are left to the reader---or see (\cite{split}).

      \bdef      \label{normal-filter}
         A set $\cF\sub\cP\cP(X)$ is a {\it normal filter} on
          $\cP(X)$ iff
    \begin{enumerate}
          \item (Filter)  $\cA,\cB\in\cF\Rightarrow\cA\cap\cB\in\cF,\:
                          \cA\sub\cB\sub\cP(X)  \a  \cA\in\cF \Rightarrow
                         \cB\in\cF,\m{ and }\emptyset\not \in\cF.$
          \item (Fineness) $\forall x\in X \: \set{a\sub X}{x\in a}
                       \in \cF.$
          \item (Normality) If $\cA_{x}\in\cF\m{ (for }x\in X)$
                        then the diagonal intersection $$\triangle_{x \in X}
			\cA_{x} =_{df}
                     \set{a\sub X}{(\forall x\in a)\,a\in\cA_{x}}\in\cF.$$
       \end{enumerate}
        If $\cF$ is a  filter on $\cP(X)$ then $I_{\cF} =_{df} 
          \mbox{$\set
          {\cA \sub \cP(X)}{\cP(X) \setminus \cA \in \cF}$}$ 
           is the dual ideal.
            $I_{\cF}^{+} =_{df}
         \set{\cS \sub \cP(X)}{\cS \not \in I_{\cF}}$ 
\edefin
\bfactnum
   Let $\cF$ be a normal filter on $\cP(X)$, $S\in I_{\cF}^{+}$ and
    $f$ a choice function on $S$ (for all $a\in S\m{, } f(a)\in a$).
    Then there is an $x\in X$ such that $\set{a\in S}{f(a)=x}\in I_{\cF}^{+}$.
\efactnum
\bremark
        It is easy to see that the conclusion of the above fact
       is actually equivalent to normality.
\eremark
\bfactnum         
       Let $\cF$ be a normal filter on $\cP(X)$.
          \begin{enumerate}
                \item    If $Y\sub X$ then 
                               the projection  of 
                             $\cF$ to     $\cP(Y),\;\gp_{X,Y}(\cF)$, is
                             a normal filter on $\cP(Y)$, where $\gp_{X,Y}(\cF)
         		=\set
                              {\gp_{X,Y}(\cA)}{\cA\in\cF}\m{ and }      
                                      \gp_{X,Y}(\cA)=\set{a\cap Y}{a \in \cA}$.
                \item     If $\cS\in I_{\cF}^{+}$ then
                                      $\cF\restricted\cS=\set{\cB\sub\cP(X)}
                               {(\exists\cA\in\cF)\:\cA\cap 
                             \cS\sub\cB}$ is a normal filter on $\cP(X)$.
         \end{enumerate}
\efactnum

\bremark
	If $Z \sub Y \sub X$ then $\gp_{Y,Z} \compose \gp_{X,Y} =
	\gp_{X,Z}$.
\eremark

\bdef  \label{club}
                  A set $\cC\sub\cP(X)$ is {\em club} (in $\cP(X)$)
                  iff $\exists f\map{X^{<\gw}}{X}$ such that
                  $\cC=\m{cl}_{f}$, where 
                  $\m{cl}_{f}=\set{a\sub X}{f\applied a^{<\gw}\sub a}$.
                  If $a \sub X$ then $cl_{f}(a) =_{df} $ the smallest set
		  containing $a$ that is closed under $f$.
\edefin
\bfactnum
               The filter generated by the club sets in $\cP(X)$,
                $\cC_{X}$, is a normal filter on $\cP(X)$.
\efactnum
\bdef          \label{stat}
                A set $\cS\sub\cP(X)$ is {\em stationary} (in $\cP(X)$)
                iff $\cS\in I_{\cC_{X}}^{+}$.  $\cS$ is {\em non-trivial} iff 
                 $X\not \in \cS.$  (Note that $\{X\}$ is stationary).
\edefin

\bremark
    If $\gk$ is regular and $\gl \geq \gk$ then $S=\set
     {a \sub \gl}{\card{a} < \gk \a a \cap \gk \in \gk}$ is stationary
    and $\cC_{\gl}\restricted S$ is the usual club filter
     on $\cP_{\gk}(\gl)$;  if $\gl = \gk$ then $\cC_{\gl}\restricted
     S$ is the usual club filter on $\gl$.  
\eremark

\bfactnum
    If $\cF$ is a normal filter on $\cP(X)$ then  $\cF$ is 
    countably complete.
\efactnum
\bfactnum
     If $\cF$ is a normal filter on $\cP(X)$ then $\cF$ contains
    the club filter $\cC_{X}$.
\efactnum

   \section{Towers of Normal Filters}

We say a set $\bp$ is a {\em tower} if there is a limit
ordinal $\gd$ (the {\em height} of $\bp$) and a function
$\cF^{\bp}\map{V_\gd}{V}$ such that for all $X \in V_{\gd}$, 
$\cF^{\bp}_X$ is a normal filter on $\cP(X)$ and for all 
$X \sub Y$ (both in $V_{\gd}$), $\cF^{\bp}_Y$ projects to $\cF^{\bp}_X$
and $\bp = \set{S \in V_{\gd}}{\exists X \in V_\gd \; S
\in I^{+}_{\cF^{\bp}_X}}$. (We often drop the superscript from
$\cF^{\bp}$.)
We define a partial order on $\bp$ by $S_1 \leq S_2$ iff
$\cup S_1 \contains \cup S_2$ and $(\forall a \in S_1) \;
a \cap (\cup S_2) \in S_2$.

In (\cite{proc}) Woodin uses the full non-stationary tower
$\Pd$ and the countable version $\Qd$.  In the above notation $\Pd$
is the tower of height $\gd$ with $\cF_X = \cC_X$ (the club
filter); $\Qd$ is the tower of height $\gd$ with $\cF_X =
\cC_X \restricted S_X$ where $S_X = \set{a \sub X}{ \card{a} \leq 
\gw}$.

\blemma \label{towerdef2}
	Assume $\bp$ is a tower of height $\gd$
          and $X,Y \in V_\gd$ with $X \sub Y$.
	Let $\gp \map{\cP(Y)}{\cP(X)}$ be the projection map $(\gp(a) =
 	a \cap X)$. Then
	\begin{enumerate}
		\item
			If $S \in I_{\cF_{Y}}^{+}$ then $\gp\applied
			S \in I_{\cF_{X}}^{+}.$
		\item
			If $S \in I_{\cF_{X}}^{+}$ then $\gp^{-1}(S)
			\in I_{\cF_{Y}}^{+}.$
	\end{enumerate}
\elemma
\bproof
	To see (1) 
 	 let $S \in I_{\cF_{Y}}^{+}$ and $C \in \cF_{X}$.  Then $\gp^{-1}
	(C) \in \cF_{Y}$ (since $\cF_{Y}$ projects to $\cF_{X}$)
  	 so there is an $a \in  \gp^{-1}(C) \cap S$ and 
	so $\gp(a) \in  C \cap \gp \applied S $. The proof of (2) is similar.
\eproof

If we force with a tower then we can form a generic ultrapower:
\blemma
	Assume $\bp$ is a tower of height $\gd$
           and $G \sub \bp$ is generic. For
	$X \in V_\gd $ let $G_{X} = \set{S \in I_{\cF_{X}^{+}}
	}{S \in G}$.  Then $G_{X}$ is
	a $V$-normal ultrafilter on $\cP(X)$ extending $\cF_{X}$.
	If $X \sub Y$  then $G_{Y}$ projects to
	$G_{X}$.
\elemma
\bproof
	Easy density arguments  show (using the above Lemma)
        that each $G_{X}$ is
	a $V$-normal ultrafilter on $\cP(X)$.  So by the definition of
	$\bp, G_{X}$ extends $\cF_{X}$.  To see projection suppose
	$X \sub Y$ are both in $V_\gd$ (and $\gp$ is the projection map).
	If $S \in G_{Y}$ then $S \leq \gp(S)$,
	 so $\gp(S)
	\in G_{X}$.  Since they are $V$-ultrafilters, $G_{Y}$ projects
         to $G_{X}$.
\eproof

So if $\bp$ is a tower of height $\gd$
 and $G \sub \bp$ is generic then we may form
the usual (direct limit) ultrapower $(M,E)$: If $f_{i}\map
{\cP(X_{i})}{V} \;(i\in \{1,2\}, f_{i} \in V,X_{i} \in V_\gd) $ then
$f_{1} \sim f_{2}$ iff for some (any) $Z \in V_\gd$ with
$X_{1} \cup X_{2} \sub Z$
$$
     \set{a\sub Z}{f_{1}(a \cap X_{1}) = f_{2}(a \cap X_{2})} \in G_{Z}.
$$
$M$ is the collection of all equivalence classes and $[f_{1}]
E[f_{2}]$ iff
$$
   \set{a\sub Z}{f_{1}(a \cap X_{1}) \in f_{2}(a \cap X_{2})} \in G_{Z}.
$$
As usual we get an elementary embedding 
$j\map{V}{(M,E)}$
and {\L}o\'{s}' Theorem: 
$j(x) = [c_{x}]$ where $c_{x} $ is the constant function $x$
with domain $\cP(Y)$ for some $Y\in V_\gd$ and
$M \models \phi([f_{1}], \dots [f_{n}])$ iff for some (any)
$Z \in V_\gd$ such that $X_{1} \cup \cdots \cup X_{n} \sub Z$
$$
   \set{a\sub Z}{\phi(f_{1}(a \cap X_{1}) \dots f_{n}(a \cap X_{n}))} \in G_{Z}.
$$
Also note that by normality for all $X \in V_\gd$,
$[f] E [\m{id}_X]$ iff there is an $x \in X$ such that
$[f] = j(x)$ ($\m{id}_X$ is the identity function with domain
$\cP(X)$).

Given any $X \in V_\gd$ we can also form the ultrapower using only $G_X$  to
get an elementary embedding $j_X\map{V}{\m{Ult}(V,G_X)}$. 
As usual, there is an elementary embedding $k\map{\m{Ult}(V,G_X)}
{(M,E)}$ defined by $k([f]_X) = [f]$.  Note that if $X$ is transitive
then $k \restricted X = \m{id}$.

\bdef
	A tower $\bp$ is {\em precipitous} if the generic ultrapower
	in $V^{\bp}$ is well-founded.
\edefin

\bdef
	Let $\bp$ be a tower of height $\gd$.
  If $A,B\sub \bp$ are antichains then
	$A \elem B$ means $(\forall p \in A) (\exists q \in B) \;p \leq q$.
	For $p,q \in \bp$ $p \sim q$ means $p \forces \quot{q \in G}$
	and $q \forces \quot{p \in G}$.  We say $p$ and $q$ are {\it disjoint}
	if for any $Z \in V_{\gd}$ with $\bigcup p , \bigcup q \sub Z$,
          $\gp_{Z,\cup p}^{-1}(p) \cap \gp_{Z,\cup q}^{-1}(q)  = \emptyset$.
	An antichain is disjoint if every pair of  elements from 
	it are disjoint.
\edefin

For a normal filter $\cF$ being able to refine every antichain in $I^{+}_
{\cF}$ to a disjoint antichain
is a strong statement (see \cite{potent}).
But for towers we have the following.

\blemma \label{disjoint}
	Assume $\bp$ is a tower of height $\gd$, $\gd$ 
 	inaccessible. If $A \sub \bp$ is a maximal antichain
	then there is a disjoint maximal antichain $B \elem
	A$.  Moreover, if $q \in B$ and $p \in A$ and $q \leq
        p$ then $q \sim p$.
\elemma
\bproof
	Let $A \sub \bp$ be a maximal antichain and $\seq{S_{\ga}:
	\ga \in \gl}$ a 1-1 listing of $A$ (so  $\gl \leq \gd$).
	It is enough to define  by induction a disjoint sequence 
	$\seq{S^{'}_{\ga}:
	\ga \in \gl}$ such that 
         $S_{\ga}^{'} \leq S_{\ga} \m{ and } S^{'}_{\ga} \sim
	 S_{\ga}$.  To define 	$S^{'}_{\gb}$
 	($\gb < \gl$)  choose $\gn < \gd$
	such that $\gn > \gb$ and 
	 for all $\ga < \gb, \cup S^{'}_{\ga} \sub V_{\gn}$
	and $\cup S_{\gb} \sub V_{\gn}$.  Let $\gp_{\ga} =
 	\gp_{V_{\gn},\cup S^{'}_{\ga}} 
	\;(\m{for }\ga < \gb)$ and $\gp_{\gb} =
 	\gp_{V_{\gn},\cup S_{\gb}} $.  Since each $S^{'}_{\ga} \sim
	S_{\ga}$ we have for each $\ga < \gb$ a set 
	$C_{\ga,\gb} \in \cF_{V_{\gn}}$ such that
	$C_{\ga,\gb} \cap \gp_{\ga}^{-1}(S^{'}_{\ga}) \cap
	\gp_{\gb}^{-1}(S_{\gb}) = \emptyset$.  
	Let $$S^{'}_{\gb} = \set{a\sub V_{\gn}}{a \in \gp_{\gb}^{-1}(S_\gb)
	\a (\forall \ga \in a \cap \gb) \:  a \in C_{\ga, \gb}
	\a \gb \in a}.$$  This clearly works.
\eproof

\vspace{.1in}
The equivalence of (1) and (2) in the following is standard (see \cite{precip}
and \cite{potent}).  
Their  equivalence with  (3) uses the above Lemma.

\blemma \label{game}
	Let $\bp$ be a tower of height $\gd$, $\gd$ inaccessible. Then the
 following are equivalent.
    \begin{enumerate}
	\item
		$\bp$ is precipitous.
	\item
		Player I does not have a winning strategy in the following
		game:  I and II alternately play elements of $\bp$ such
		that $p_{0} \geq p_{1} \geq \cdots \geq p_{n} \cdots$
		and II wins iff $\exists a \sub \bigcup_{n\in \gw}
		(\cup p_{n})$
		such that $(\forall n)\; a \cap (\cup p_{n}) \in p_{n}$.

	\item
		If $p \in \bp$ and $A_{n} \sub \bp$ are maximal
		antichains ($n\in \gw$) then 
$$
     \set{a \sub V_{\gd}}
		{a \cap \cup p \in p \a (\forall n)(\exists q \in A_{n})\;
		a \cap \cup q \in q} \not =  \emptyset.
$$
		
   \end{enumerate}
\elemma
\bproof
	(1)$\!\implies\!$(2).  Assume that (2) fails.  Let $\gs$ be a winning 
    	strategy for I.  Define a tree $T$:  $\seq{(p_0,a_0), \dots (
   (p_n,a_n)} \in T$ iff
\begin{enumerate}
\item $\seq{p_0, \dots p_n} $ is according to $\gs$.
\item $(\forall i \in n)\; a_i  \in p_i$.
\item $(\forall i < j \leq n)\; a_j \cap (\cup p_i) = a_i$.
\end{enumerate}
Since $\gs$ is a winning strategy for I, $T$ is well-founded.
Let $G \sub \bp$ be generic with $\gs(\emptyset) \in G$.  So there
exists (in $V[G]$) a sequence $\seq{p_0,p_1,\dots}$ according to $\gs$
such that every $p_i \in G$.  Let $j\map{V}{M} \sub V[G]$ be the
generic embedding.  If  $M$ is well founded
then $j(T) $ is well founded in $V[G]$.  But
$$ \seq{ (j(p_0),j\applied \cup p_0), \dots 
 (j(p_n),j\applied \cup p_n), \dots }
$$
is an infinite descending chain in $j(T)$.

	(2)$\implies$(3).
	Let $p, A_{n}$ witness the failure of (3).  A winning strategy for
	I is to let $\gs(\emptyset) = p$ and at the $n$th move play something
	below an element of $A_{n}$ (and below II's last move).
	
	(3)$\implies $(1).
Assume (1) fails.  So there is a $p\in \bp$ and names $\gt_n$ such that
$(\forall n \in \gw) \; p \forces \quot{\gt_n,\gt_{n+1} \in \m{Ord} \a
\gt_n > \gt_{n+1}}$.  Let $A_{-1}$ be any disjoint maximal antichain with
$p \in A_{-1}$.  Inductively construct disjoint maximal antichains $A_n$
and functions $T_n$ such that
\begin{enumerate}
\item $A_{-1} \succ A_0 \succ A_1 \dots$
\item $(\forall n \in \gw) \;  \forall q \in A_n$, if $q \leq p$ then
$T_n(q)\map{q}{\m{Ord}}$ such that $q \forces \gt_n = [T_n(q)]$.
\item  Suppose $p \geq q_1 \geq q_2$ and $q_1 \in A_{n_1}$, 
$q_2 \in A_{n_2}$, and $n_2 > n_1$.  Then ($\forall a \in q_2$)
$T_{n_2}(q_2)(a) < T_{n_1}(q_1)(a \cap (\cup q_1)$.
\end{enumerate}

But now $p$, $\seq{A_n: n \in \gw}$ witnesses the set defined in (3)
is empty:  Suppose $a \sub V_\gd$ is in this set.  Then $a \cap
(\cup p) \in p $ and $(\forall n)\exists q_n \in A_n$ such that
$a \cap (\cup q_n) \in q_n$.  Our construction gives that
$p \geq q_0 \geq q_1 \dots$, and so
$T_{0}(q_0)(a \cap (\cup q_0)) >
T_{1}(q_1)(a \cap (\cup q_1)) > \dots$
\eproof

\bremarkss
The equivalence of (1) and (2) does not use that $\gd$ is inaccessible.

     We can also add to the above list: 

		(4)  If $p \in \bp$ and $A_{n} \sub \bp$ are maximal
		antichains ($n\in \gw$) then 
$$
     \set{a \sub V_{\gd}}
		{p \in a \a 
                      a \cap \cup p \in p \a (\forall n)(\exists q \in A_{n}
		\cap a )\;
		a \cap \cup q \in q} \not =  \emptyset.
$$
If $I(\gd)$ is the above set then $I(\gd) \not = \emptyset$ is
also equivalent to $\exists \gk < \gd$ such that $I(\gk) \not
= \emptyset$ (since cof$(\gd) > \gw$).  This form of precipitous
is used in \cite{proc}.
Also note that if for all $p \in \bp$ and all maximal antichains $A_n \sub
\bp$ ($n \in \gw$) $\exists \gk < \gd$ such that $I(\gk) \in I_{
\cF_{V_\gk}}^{+}$ then the generic ultrapower is
closed under $\gw$-sequences in $V[G]$.  In fact, this is
equivalent to being closed under $\gw$-sequences.  So, for example, if
$\bp$ is the tower of height $\gd$ where $\cF_X = \cC_X \restricted
S_X$ and $S_X = \set{a\sub X}{\card{a} < \aleph_\gw}$ then $\bp$ is
precipitous (if $\gd$ is Woodin) but the generic ultrapower is 
not closed under $\gw$-sequences.
\eremarkss

   \section{Well-founded extensions of filters}

In  this section we show that any normal filter is part
of a tower of arbitrarily large height.  So (assuming large cardinals) 
every normal filter can be generically extended to a 
well-founded $V$-ultrafilter --- although
the filter itself may not be precipitous. We also prove a similar result
for any tower on $V_\gd$, assuming that $\gd$ is inaccessible.

 We use the following
 Lemma---Foreman proved this when when $\cF$ is an ultrafilter.
\blemma \label{fwf}
   Assume $\cF$ is a normal filter on $\cP(Y)$.  Assume 
that $\card{X} \geq 2^{2^{\card{Y}}}$ and $Y \sub X$.
     Then there is a stationary set $S$
in $\cP(X)$ such that the club filter on
$\cP(X)$ restricted to $S$ projects to $\cF$.
\elemma
\bproof
   Let $\gp$ be the projection map from $\cP(X)$ to $\cP(Y)$.\
Let $\seq{C_{x}: x \in X}$ be a listing of all
the elements in $\cF$.  Let 
$$ S = \set{a \sub X}{\forall x \in a \:  (a \in \gp^{-1}
(C_{x}))}.$$
We need to see that $S$ is stationary and that
$\gp(\cC_{X} \restricted S) = \cF$.  For this,
it is enough to show that $\forall f \map
{X^{<\gw}}{X}$, $\gp(S \cap \m{cl}_{f}) \in \cF$.
Fix such an $f$.
Let $\bar{f}\map{Y^{<\gw}}{Y}$ be such that if $a \sub Y$ is closed
 under $\bar{f}$ then $\m{cl}_{f}(a) \cap Y = a$ (set $\bar{f}(y_{1} \dots
y_{\seq{i,j}}) = $ the $j$'th element of $Y$ in
cl$_{f}(y_{1} \dots y_{i})$, where $\seq{\cdot,\cdot}$ is some simple
pairing function on $\gw$).
\bclaim
  Suppose that $g\map{Y^{<\gw}}{\cP_{\gw_{1}}(\cF)}$.  Then
$$
\set{a \sub Y}{ \forall \gt \in a ^{<\gw} \: \forall C \in g(\gt)
(a \in C)} \in \cF.$$
\eclaim
\bproofclaim
Let $h\map{Y^{<\gw}}{Y}$ be a bijection.
For $y \in Y$ let $D_{y} = \bigcap g(h^{-1}(y))$
(so $D_y \in \cF$).  Then
cl$_{h} \cap \triangle_{y \in Y} D_{y} \sub
\set{a \sub Y}{ \forall \gt \in a ^{<\gw} \: \forall C \in g(\gt)
(a \in C)} $,  so this set is in $\cF$.
\eproofclaim
Now given $\gt \in Y^{<\gw}$ let $g(\gt)=\set
{C_{x}}{x \in \m{cl}_{f}(\gt)}$.
So 
$$
C=\set{a \sub Y}{a \in \m{cl}_{\bar{f}} \a \forall \gt \in
a^{<\gw} \: \forall C \in g(\gt) \: (a \in C)} \in \cF.$$
But if $a \in C$ then cl$_{f}(a) \in S \cap \m{cl}_{f}$ and 
$\m{cl}_{f}(a) \cap Y = a$.
Therefore $\gp(S\cap \m{cl}_{f}) \in \cF$.
\eproof
As a corollary to this and a Theorem of Woodin we get the following.
\bcor
Assume $\cF$ is a normal filter on $\cP(X)$ and there is a
Woodin cardinal $> \card{X}$.  Then in some generic extension of $V$, there
is a $V$-ultrafilter extending $\cF$ with
well-founded ultrapower.
\ecor
\bproof
  We use the result of Woodin (\cite{cabal}) 
that if $\gd$ is Woodin and $G \sub \Pd\
$ is generic then the direct limit ultrapower is
well-founded (and so the ultrapower using any
measure from $G$ is well-founded).
Let $S$ be a stationary set on some $\cP(Y)$ ($Y \in V_{\gd}$) such that
$\cC_{Y} \restricted S$ projects to $\cF$
(we may assume that $X$ is an ordinal, so $X \in V_{\gd}$).
  Let $G \sub \Pd\
$ be generic with $S \in G$.  Then
$$ \set{S^{'}}{S^{'}  \in G \a S^{'} \m{ is stationary in } \cP(X)}$$
 is a $V$-ultrafilter extending
$\cF$ with well-founded ultrapower.
\eproof

\bthm
Assume $\bp$ is a tower on $V_\gd$, $\gd$ inaccessible.
 Then there is a stationary set $S$ in $\cP(V_\gd)$ such that
for all $X \in V_\gd$, $\cC_{V_\gd} \restricted S$ projects
to $\cF_X$.
\ethm

\bproof
Let $S = \set{a \sub V_\gd}{(\forall X \in a) (\forall C \in
\cF_X \cap a) \; a \cap X \in C}$.  It is enough to see that for every 
$X \in V_\gd$ and every $f\map{V_{\gd}^{< \gw}}{V_\gd}$,
$\gp_{V_\gd,X}(\m{cl}_f \cap S) \in \cF_X$.

Fix $X$ and $f$.  Since $\gd$ is inaccessible, 
 $\exists \gb < \gd$ such that $V_\gb$ is closed under $f$
(and $X \in V_\gb$).  But then 
$$\overline{S} = 
\set{a \sub V_\gb}{a \in \m{cl}_f \a (\forall Y \in a) (\forall C
\in \cF_Y \cap a) \; a \cap Y \in C} \in \cF_{V_\gb}$$
(Since $\cF_{V_\gb}$ projects to $\cF_Y$ for all $Y \in
V_\gb$.)  So we are done: the projection of $\m{cl}_f \cap
S$ to $V_\gb$ contains $\overline{S}$ and the projection
of $\overline{S}$ to $X$ is in $\cF_X$.
\eproof

   \section{Examples of non-precipitous towers}

In this section we give examples (which were suggested by Woodin)
 of non-precipitous towers
(assuming the existence of a supercompact).
These examples use Lemma \ref{notin} below, which says that
under certain conditions  (precipitousness and moving its height)
towers are not in their ultrapowers.
We do not know if these conditions are necessary nor if the
supercompact is needed for these examples.

The proof of the following Lemma is based on a proof of the fact that
ultrafilters are not in their ultrapowers.

\blemma \label{notin}
Assume $\bp$ is a precipitous tower of height $\gd$, $\gd$ 
inaccessible, and $V^{\bp}
\models i_G(\gd) > \gd$.  Then $\bp$ is not in its generic ultrapower.
\elemma

\bproof
Assume the Lemma fails.  Let $G \sub \bp$ be generic and
$j\map{V}{M} \sub V[G]$ the generic embedding with $j(\gd) >
\gd$ and $\bp \in M$.

Note that $G$ is also generic over $L(\bp)$ and that $V_\gd 
\in L(\bp)$, so $\exists p \in \bp$ such that 
$L(\bp) \models \quot{p \forces i_G(\gd) \geq
\check{j(\gd)}}$  (since $G$ witnesses it).
Let $[d] = \gd$ and $[p] = \bp$.  We may assume $\dom(d) = \dom(p)
= V_{\ga + 1}$ $(\ga < \gd)$ and (since $j(\gd) > \gd$)
($\forall a \sub V_\ga$) $d(a) < \gd$ and $p(a) \sub V_{d(a)}$
is a tower in $L(p(a))$ of height $d(a)$.
But then $L(p(a))^{p(a)} \models \quot{i_G(d(a)) < \check{\gd}}$
(since $\gd$ is inaccessible).  So in $M$, $L(\bp)^{\bp} \models
i_G(\gd) < \check{j(\gd)}$.  Contradiction.
\eproof

Now assume that $\gk$ is supercompact and $\gd > \gk$ is inaccessible.
We will define a tower of height $\gd $ that is not precipitous.
Let $$A_0 = \set{\gu}{(\exists X \in V_\gd)\; \gu \m{ is a supercompactness
measure on } \cP_{\gk}(X)}$$  For $\gu \in A_0$ let supp$(\gu) = $ the
unique $X$ such that $\gu$ is a supercompact measure on $\cP_{\gk}(X)$.
Inductively define $A_\gl$: for limit $\gl$, $A_\gl = \bigcap_{
\ga < \gl} A_\ga$; given $A_\gl$, let $A_{\gl + 1} = \set{
\gu \in A_\gl}{(\forall Y \in V_\gd) \m{ if } Y \contains \m{supp}(\gu)
\m{ then } (\exists
\gn \in A_\gl) \; \m{supp}(\gn) = Y \a \gn \m{ projects to } \gu}$.
Let $A = \bigcap A_\ga$.  Note that $A$ is non-empty:  the projections
of any supercompactness measure on $\cP_{\gk}(V_\gd)$ are in $A$.
By construction $
(\forall \gu \in A) (\forall Y  \in V_\gd)
\m{ if } Y \contains \m{supp}(\gu) \m{ then }
 (\exists
\gn \in A) \; \m{supp}(\gn) = Y \a \gn \m{ projects to } \gu$.
Also note that the measures in $A$ are closed under projection.
Given $X \in V_\gd$, let $\cF_X$ be the filter on $\cP(X)$ generated by
$\set{\cC \sub \cP_\gk(X)}{(\forall \gu \in A) \m{ if supp}(\gu) = X \m{ then }
\cC \in \gu}$.  These are normal filters that project to one another;
let $\bp$ be the associated tower.  Assume that $\bp$ is precipitous.
Note that $i_G(\gk) \geq \gd$ and so $i_G(\gd) > \gd$.  We
will get a contradiction by showing
that $\bp$ is in the generic ultrapower. Let $G \sub \bp$ be generic
and $j\map{V}{M} \sub V[G]$ the generic embedding.

\bclaim  $(\forall X \in V_\gd)\; \exists \gu \in A$ such that the generic
ultrafilter $G_X = \gu$.
\eclaim
\bproofclaim
Fix $X \in V_\gd$.  We may assume $X = V_\ga$ for some $\ga$.  Fix
$S \in \bp$.  By extending $S$ if necessary we may assume that
$\cup S = V_\gb \contains V_{\ga +2}$ and $ S \sub \cP_\gk(V_\gb)$.
Since $ S \in \bp$ there is a $\gu \in A$ such that supp$(\gu) =
V_\gb$ with $S \in \gu$.  By Lemma \ref{fwf} there is a 
stationary $S^{'}$ in $\cP(V_\gb)$ such that $\cC_{V_\gb} 
\restricted S^{'}$ projects to $\gu \restricted V_\ga$.  The proof
of Lemma \ref{fwf} shows that $S^{'} \in \gu$.  Hence $S \cap
S^{'} \in \bp$ and $S \cap S^{'} \forces G_{V_\ga} = \gu \restricted
V_\ga$.
\eproofclaim
\bclaim  $M_\gd = V_\gd$
\eclaim
\bproofclaim
It is always the case that $V_\gd \sub M_\gd$.  So let $\cA \in M_\gd$.
Since $M$ is the direct limit of the $j_X\map{V}{\m{Ult}(V,G_X)}$,
there is an $\ga < \gd$
 such that $\cA \in M_\ga$ and an $\overline{\cA}
\in \m{Ult}(V,G_{V_\ga})$ such that $k(\overline{\cA}) = \cA$
(where $k$ is the canonical map from $\m{Ult}(V,G_{V_\ga})$ into
$M$).  But $k \restricted V_\ga = \m{id}$, so $k(\overline{\cA}) =
\overline{\cA}$.  Hence $\cA \in \m{Ult}(V,G_{V_\ga}) \sub V$.
\eproofclaim
But our construction of $\bp$ is absolute to $L(V_\gd)$, so 
$\bp \in L(V_\gd) \sub M$.  Contradiction, so $\bp$ is not
precipitous.
 
   \section{Large Cardinals}

In this section we will describe the large cardinal we use.
A cardinal $\gd$ is $\gl$-{\em supercompact} if there is an elementary
embedding $j\map{V}{M}$ such that c.p.($j)=\gd$ and
$M^{\gl} \sub M$.  
For $A$ a set of ordinals we say $\gd$ is $[\gl] \: A-${\em superstrong}
if there is an elementary embedding $j\map{V}{M}$ such that c.p.($j)=\gd$,
$j(A) \cap j(\gd) = A \cap j(\gd)$
and $j \applied \gl \in M$.

\bthm   \label{thm.large.cardinal}
	Assume $\gd$ is $\card{V_{\gd + \gw + 2}}-$supercompact.  Then
	for all $A \sub \gd$ there are stationary many $\gk < \gd$ such
	that $\gk$ is $[\card{V_{\gk + \gw}}]\; A-$superstrong.
\ethm

\bremarkss
       1.  Actually, all we need (in Theorem \ref{thm.semiproper})
           is some small amount of ``$[\card{V_{\gk + \gw}}]\; A-$strong''
           for a certain $A$.

       2.  The proof of this Theorem follows the proof of ``if $\gd$ is
	$2^{\gd}-$supercompact then $\gd$ is a Woodin cardinal''
	(see \cite{ms}).  
	
       3.  In fact, if we let $\gl(\ga) = \card{
	V_{\ga+\gw}}$ (or other simple functions like $\gl(\ga)=
	2^{\ga}$ or $\gl(\ga) = \ga$) then this same method of proof
	gives the following:
   \begin{list}{}{\setlength{\leftmargin}{.5in}}
	\item[a.]
		If $\gd$ is supercompact then there are
		stationary many $\gk < \gd$ such that $\gk$ is
		$[\gl(\gk)] $ superstrong. 
	\item[b.]
		If $\gd$ is $[\gl(\gd)]$ superstrong then there are
		stationary many $\gk < \gd$ such that $\gk$ is
		$[\gl(\gk)] $ Shelah and $\gd$ is $[\gl(\gd)] $ Shelah
		(where $\gk$ is $[\gl(\gk)] $ Shelah has the obvious
		definition).
	\item[c.]
		If $\gd$ is $[\gl(\gd)]$ Shelah then there are
		stationary many $\gk < \gd$ such that $\gk$ is
		$[\gl(\gk)] $ Woodin (and $\gd$ is $[\gl(\gd)] $ Woodin).
   \end{list}
 	The case $\gl(\ga) = \ga$ is what is proved in \cite{ms}.

       4.  We can strengthen these results  by
          requiring that $M^{\gl} \sub M$
	   rather than just $j \applied \gl \in M$.
\eremarkss
\begin{pf*}{Proof of Theorem \ref{thm.large.cardinal}} 
	 Assume the
	Theorem fails.  So there is an $A \sub \gd$ and a club $C \sub
	\gd$ such that if $\gk \in C$ then $\gk$ is not
	$[\card{V_{\gk + \gw}}]\;  A-$superstrong.
	Let $j\map{V}{M}$ with cp($j) = \gd$ and $M^{V_{\gd + \gw +2}}
	\sub M$. So $$ (*) \ \ \ M \models \quot{\gd \m{ is not }
	[\card{V_{\gd + \gw}}]\; j(A)-\m{superstrong}}$$
	Let $E$ be the sequence of
        measures  derived from j with support
	$ S = j(\gd) \cup \{j \applied V_{\gd + \gw} \}$.
       ($E = \seq{\gu_{\gt}: \gt \in S^{<\gw}}$ where $\gu_{\gt}
       (X) = 1 $ iff $\gt \in j(X)$.)
       We can form the (direct limit)
       ultrapower using $E$:  $(f,\gt) \sim (g,\gs)$ (for $
       \gt,\gs \in S^{<\gw}) $ iff $j(f)(\gt) = j(g)(\gs)$ and
       $[f,\gt] E [g,\gs] $ iff $j(f)(\gt) \in j(g)(\gs)$
      (so the ultrapower is well-founded).
       Let $i_{E}\map{V}{\m{Ult}(V,E)}$ be the canonical embedding.
       We will show that $E \in M$ and use this to contradict $(*)$.
         \bclaim
	   Define $id^{*}(\seq{a}) = a$.  For any $\ga < j(\gd)$,
	   $[id^{*},\seq{\ga}] = \ga$.
         \eclaim
	 \bproofclaim  Easy, by induction on $\ga$.
	 \eproofclaim
	 \bclaim
		$i_{E}(A) = j(A)$.
         \eclaim
	 \bproofclaim
     Using the above claim it is easy to see that 
	     $i_{E}(\gd)  = j(\gd)$.
	       So $\ga \in i_{E}(A)$ iff $[id^{*},\seq{\ga}]
	     \in i_{E}(A)$ iff $\ga \in j(A)$.
         \eproofclaim
	 \bclaim
	      $E \in M$.
          \eclaim
	  \bproofclaim
Since for every $\gt \in S^{<\gw}$, $\gu_\gt(V_{\gd + \gw +1}^{<\gw}) = 1$,
	      $E$ is defined from its support and $j \restricted 
	      \cP(V_{\gd + \gw +1}^{<\gw})$. Since 
$M$ is closed under $\card{V_{\gd + \gw +2}}$ sequences, $E$ is in $M$.
           \eproofclaim

	   So we can form the ultrapower of $M$ by $E$ and get an elementary
	   embedding $i_{E}^{M}\map{M}{\m{Ult}(M,E)}$.  This ultrapower
	   is well-founded:
	   \bclaim
	       $[f,\gt]_{M} \in [g,\gs]_{M}$ (for $f,g \in M$) iff $j(f)(\gt)
	       \in j(g)(\gs)$.
            \eclaim
	    \bproofclaim
	       $[f,\gt]_{M} \in [g,\gs]_{M}$ iff
               (by definition) $ \set{ a_{1} \concat a_{2}}{
	       f(a_{1}) \in f(a_{2})} \in \gu_{\gt \concat \gs}$ iff
	       $j(f)(\gt) \in  j(g)(\gs)$.
            \eproofclaim
	    Finally, the following claim contradicts $(*)$.
	    \bclaim
	      $i^{M}_{E}$ has critical point $\gd$,
	      $i^{M}_{E}(\gd) = j(\gd)$,
	      $i^{M}_{E}(A) = j(A)$ and
	      $i^{M}_{E} {\applied} V_{\gd + \gw} \in \m{Ult}(M,E)$
             \eclaim
              \bproofclaim
                 $M$ and $V$ have the same $V_{\gd +\gw +2}$ so they
                have the same functions from
		$V_{\gd + \gw +1}^{<\gw}$ into
		$V_{\gd+ 20}$.  Hence
                 c.p.($ i_{E}^{M}) = \gd$,
            $ i_{E}^{M}(\gd) = i_{E}(\gd)$ and
            $ i_{E}^{M}(A) = i_{E}(A)$.  
              Finally, it is easy to see that $[id^{*},\seq{j \applied
             V_{\gd + \gw}}]_{M} = i^{M}_{E} {\applied} V_{\gd + \gw}$.
               \end{pf*} 
          \renewcommand{\qed}{}      \end{pf*} 

   \section{Proof of Precipitousness}

\bthm \label{thm.precip}
     Assume that $\bp$ is a tower of height $\gd$  where $\gd$
is $\card{V_{\gd + \gw +2}}$--supercompact and for all $X \in V_{\gd}$,
$\cF^{\bp}_{X} = \cC_{X} \restricted S_{X}$ for some stationary set $S_{X}$.
Then $\bp$ is precipitous.
\ethm

\bdef
	A set $b$ {\em end extends} $a$ if for all $x \in a, a \cap
	x = b \cap x$.
\edefin

\bdef
	Let $\bp$ be a tower of height $\gd$.
	For $\gk < \gd$,
	$A \sub
	\bp \cap V_{\gk}$,
	$\gl $ any ordinal $> \gk^{+}$ and $s$ any set in $V_{\gk^+}$
       we define
	$C_{s}(\gk,\gl,A)$ 
 to be the set of all $a \elem V_{\gk + \gw}$ such that:
   given any $a^{*}$ (with $a^{*}  \elem a, \card{a^{*}} < \gk, 
   a \m{ end extends } a^{*} \cap V_{\gk},
               A \in a^{*}$)  there exists $
              b$ such that
      \begin{enumerate}
        \item   $b \elem V_{\gl}$ with   
                 $\cF_{V_{\gk+\gw}} \in b$, $s \in b$ and
                   $a^{*} \sub b$.
	\item
		$\forall C \in \cF_{V_{\gk+\gw}}\; (C \in b \implies 
		b \cap V_{\gk + \gw} \in C)$.
	\item
		$b$ end extends $a^{*}  \cap V_{\gk}$.
	\item
		$(\exists x \in b \cap A)\; b \cap \cup x \in x$.
   \end{enumerate}
\edefin

\bthm \label{thm.semiproper}
	Assume $\bp$ is a tower of height $\gd$
	 and $\gd$ is $\card{V_{\gd+\gw+2}} $--supercompact.
	Then there are stationary many inaccessible
	cardinals $\gk < \gd$ such that
	for any ordinal $\gl > \gk^{+}$, any $s \in V_{\gk^+}$
       and   any maximal antichain $A \sub
	\bp \cap V_{\gk}$, the set $C_{s}(\gk,\gl,A)$  
is in $\cF_{V_{\gk+\gw}}$.
\ethm

\bproof
	Assume the Theorem fails.  So there is a club $C \sub \gd$ such
	that if $\gk$ is inaccessible and in $C$ then there is a $\gl >
     \gk^{+}$, $s \in V_{\gk^+}$
	and a maximal antichain $A \sub \bp \cap V_{\gk}$ 
	such that $C_{s}(\gk,\gl,A) \not \in \cF_{V_{\gk + \gw}}$.
	By Theorem \ref{thm.large.cardinal} 
	there is a $\gk \in C$ and an elementary
	embedding $j\map{V}{M}$ with critical point $\gk$
	such that $j \applied V_{\gk+\gw} \in M$ , $V_{\gk +\gw + \gw}
	\sub M$ and $j(\bp) \cap V_{\gk+\gw+\gw} = \bp \cap
	V_{\gk+\gw+\gw}$.  Fix a  $\gl$, $s$ and $A$ for
        this $\gk$,  so
	$$p = \set{a \elem V_{\gk+\gw}}{a \not \in C_{s}(\gk,\gl,A)} \in
	I^{+}_{\cF_{V_{\gk+\gw}}}.$$  So $p \in M$ and $p \in
	\bp \cap V_{\gk +\gw +\gw}$ and therefore $p \in j(\bp \cap
	V_{\gk})$.  So there is a $q \in j(A)$ and $r \in j(\bp \cap
	V_{\gk}) $ such that $r \leq p,q$.  
	
	Let $b \elem V^{M}_{j(\gl)}$, $b \in M$, such that
    \begin{enumerate}
	\item
		$j(\gk),p,q,r,j(\cF)_{V_{j(\gk)+\gw}},j(s),j\restricted
		V_{\gk+\gw} \in b$
	\item 
		$\forall C \in j(\cF)_{V_{j(\gk)+\gw}}(C \in b \implies
		b \cap j(V_{\gk+ \gw}) \in C)$.
	\item
		$b \cap \cup r \in r$.
    \end{enumerate}
	Since $r \leq p, a = _{df} b \cap V_{\gk + \gw} \in p$.  Therefore
	$\exists a^{*} \elem a \: (\card{a^{*}} < \gk \a
        a \m{ end extends } a^{*} \cap V_{\gk} \a
        A \in a^{*})
	$ such
	that no $  b$ satisfies conditions
	1-4 in the definition of $C_{s}(\gk,\gl,A)$. 
	Fix  a witness $a^{*}$.  So in $M$ the same must be true of
	$j(a^{*})$ and conditions 1-4 in $C_{j(s)}(j(\gk),j(\gl),j(A))$.
	Now use $b$ from
	above 
	 to get a contradiction in $M$. 
        Note that $j(a^{*}) = j \applied a^{*}$ (since $\card{a^{*}} < \gk$)
	so $j(a^{*}) \sub b$ (since $j \restricted V_{\gk + \gw} \in b$)
	and $b$ end extends $j(a^{*}) \cap j(V_{\gk})$ (since
	$b \cap V_{\gk + \gw} =a$ which end extends $a^{*} \cap V_{\gk})$.
	Finally, we get condition (4) since $q \in b \cap j(A)$ and
	$b \cap \cup q \in q $ (since $r \leq q$).
\eproof

\bremark
    We can get by with much weaker assumptions on $\gd$
if we drop condition (2) from the definition of 
$C_{s}(\gk,\gl,A)$.  For example, using the notation from the above 
proof, if $M^{<\gk} \sub M$ and the set $r$ is stationary
in $V$, then we also reach a contradiction:  We find $b
\elem V_{j(\gl)}^{M}$ with $b \in V$ as above (except we
drop condition (2)).  Then $\m{Hull}^{M}(j(a^{*}) \cup \{q\} 
\cup \{j(s)\} \cup \{ \cF_{V_{\gk + \gw}} \}
\cup (b \cap \cup q) ) \in M $ plays the role of $b$.

So if we let $C_{s}^{*}(\gk,\gl,A)$ be defined like $C_{s}(\gk,\gl,A)$ but we
drop condition (2),  then if $\gd$ is Woodin and $
\seq{A_\ga : \ga \in \gd}$ is a sequence of
maximal antichains in $\bp$, then there are stationary many
inaccessibles $\gk < \gd$ such that $(\forall \ga < \gk)
(\forall \gl > \gk^{+}) (\forall s \in V_{\gk^+})
 \; C_{s}^{*}(\gk,\gl,A_\ga \cap V_\gk) \in
\cF_{V_{\gk + \gw}}$.  If the tower has a ``simple'' definition
(for example $\Pd$ or $\Qd$) then condition (2)
in $C(\gk,\gl,A)$ holds automatically, so the proof of precipitousness
below goes through.
Note that the above Theorem does not need the assumption that
$\cF_X = \cC_X \restricted S_X$.
\eremark

\ \ \

\noindent {\it Proof of Theorem} \ref{thm.precip}.
	Let $\bp$ be a tower of height $\gd$,  where 
	 $\gd$ is $\card
	{V_{\gd +\gw + 2}}$--supercompact.  
  Assume ($\forall X \in V_\gd$)
 $\cF^{\bp}_{X} = \cC_X \restricted S_{X}$ for some stationary set $S_{X}$.
	We will verify condition (3) of Lemma  \ref{game}. So let $p \in
	\bp$ and $A_{n} \sub \bp$ be maximal antichains ($n \geq 1$).
	Since  there are club many $\gk < \gd$ such that
	$\forall n \;
	 (A_{n} \cap V_{\gk})$ is a m.a.c. in $\;\bp \cap V_{\gk}$, by Theorem 
	 \ref{thm.semiproper}
	 there is an inaccessible $\gk < \gd$ such that
	 $p \in V_{\gk}$ and (letting $B_{n} = A_{n} \cap
	 V_{\gk}$ and $s = S_{V_{\gk + \gw}})\;\;\; \forall n \geq~1,
	 \;\; C_{s}(\gk,\gd,B_{n}) \in \cF_{V_{\gk+\gw}}$.

	 Let $\gn >> \gd$ (say $\gn$ is strong limit, $\m{cf}(\gn)
	 > \gd$).  Choose $a_{0} \elem V_{\gn}$ such that
	 $p,\gk,\gd,\bp,s,A_{1}, \dots \in a_{0}$ and
	 $a_{0} \cap \cup p \in p$ and $a_{0} \cap V_{\gk + \gw}
	 \in C_{s}(\gk,\gd,B_{1})$.  Let $\gee_{0} < \gk$ be a limit
	ordinal in $a_{0}$ such that $p \in V_{\gee_{0}}$.  Assume inductively
	that we have defined $a_{0}, \dots a_{n}$ and $
	\gee_{0} < \cdots < \gee_{n} < \gk$ such that $a_{n} \elem V_{\gn}$ and 
	$p,\gk,\gd,\bp,s,\gee_0 \dots \gee_n,
         A_{1}, \dots \in a_{n}$ and for all $i \leq
	n$ we have (letting $B_{0} = \{p\}$) $\exists x_{i} \in
	a_{i} \cap B_{i}\; (a_{i} \cap \cup x_{i} \in x_{i} \a
	x_{i} \in V_{\gee_{i}})$ and if $i < n$ then
	$a_{i} \cap V_{\gee_{i}}=
	a_{i+1} \cap V_{\gee_{i}}$.  Also, $a_{n} \cap V_{\gk+\gw}
	\in C_{s}(\gk,\gd,B_{n+1})$.  

	If we can keep going with this construction then we are
	done since $\bigcup_{n \in \gw}(a_{n} \cap V_{\gee_{n}})$
	witnesses that the set in part (3) of  Lemma \ref{game}
	is non-empty.

	To define $a_{n+1}$ let $a_{n}^{*} \elem a_{n}$ with
	$\card{a_{n}^{*}} < \gk$, and
          	   $p,\gk,\gd,\bp,s,
	\gee_0, \dots \gee_{n},A_{1}, \dots 
	\in a_{n}^{*}$, and
         $a_{n}$ end extends $a_{n}^{*}
	\cap V_{\gk}$,
         and $ a_{n} \cap V_{\gee_{n}} \sub a_{n}^{*}$.
	Since $a_{n} \cap V_{\gk +\gw} 
	\in C_{s}(\gk,\gd,B_{n+1})$ and $a^{*}_{n} \cap V_{\gk + \gw}$ has
	the required properties, there exists a $
	 b$ such that 
   \begin{enumerate}
        \item   $b \elem V_{\gd}$ with   
                 $\cF_{V_{\gk+\gw}} \in b$, $s \in b$ and
                   $a_{n}^{*} \cap V_{\gk+\gw} \sub b$.
	\item
		$\forall C \in \cF_{V_{\gk+\gw}}\; (C \in b \implies 
		b \cap V_{\gk + \gw} \in C)$.
	\item
		$b$ end extends $a_{n}^{*}  \cap V_{\gk}$.
	\item
		$(\exists x_{n+1} \in b \cap B_{n+1})\; b \cap \cup 
		x_{n+1} \in x_{n+1}$.
   \end{enumerate}

	  Let $a_{n+1} = \set{f(\gt)}{f \in a_{n}^{*} \a
	\gt \in b \cap V_{\gk + \gw}}$. 
	We verify the inductive assumptions for $a_{n+1}$.
   \bclaim
	$a^{*}_{n} \sub a_{n+1} \elem V_{\gn}$ and $a_{n+1} \cap 
	V_{\gk +\gw} = b \cap V_{\gk + \gw} $.
   \eclaim
   \bproofclaim
	Clearly $a^{*}_{n} \sub a_{n+1}$.  To see that $a_{n+1} \elem
	V_{\gn}$ suppose $V_{\gn} \models \quot{\exists x \phi(
	x,f_{1}(\gt_{1}), \dots f_{n}(\gt_{n}))}$. Since $\gn$ is large
   $$
	V_{\gn} \models \exists g \forall x_{1}, \dots
	x_{n} \in V_{\gk + \gw} [\exists x \phi(x,f_{1}(x_{1})
	, \dots
	) \implies \phi(g(\seq{x_{1}, \dots x_{n}}),f_{1}(x_{1}),
	 \dots  )]
   $$
	So there is such a $g$ in $a_{n}^{*}$ and hence $a_{n+1} 
	\elem V_{\gn}$.

	To see $a_{n+1} \cap V_{\gk + \gw} \sub  b \cap V_{\gk + \gw}
	$ suppose that $f(\gt) \in V_{\gk + \gw}$.  Since $\gt \in
	V_{\gk + \gw}$ we may assume $f\map{V_{\gk + n}}{V_{\gk +n}}$
	for some $n \in \gw$.  But $a^{*}_{n} \cap V_{\gk + \gw}
	\sub b$ so $ f \in b$ and hence $ f(\gt) \in b$.
	The other inclusion is clear.
   \eproofclaim

	Now we need to check three conditions and to define $\gee_{n+1}$:
    \begin{enumerate}
	\item
		$a_{n} \cap V_{\gee_{n}} = a_{n+1} \cap V_{\gee_{n}}$:
		This holds since 
		$a_{n} \cap V_{\gee_{n}} = a_{n}^{*}
		\cap V_{\gee_{n}} = b \cap V_{\gee_{n}}$ (since $\gee_{n}
		\in a_{n}^{*} \cap V_{\gk}$)
		and  $b \cap V_
		{\gee_{n}} = a_{n+1} \cap V_{\gee_{n}}$.
	\item
		$(\exists x_{n+1} \in a_{n+1}\cap B_{n+1}) \:
	 	a_{n+1} \cap \cup x_{n+1}
		\in x_{n+1}$:
		This holds since there is such an $x_{n+1}$ in
		$b$ and $b \cap V_{\gk + \gw}
		= a_{n+1} \cap V_{\gk + \gw} $.
		Let $\gee_{n+1}< \gk$ be a limit ordinal in $a_{n+1}$
		such that $a_{n+1} \in V_{\gee_{n+1}}$ and
		$\gee_{n+1} > \gee_{n}$.
	\item
             $a_{n + 1} \cap V_{\gk + \gw} \in C_{s}(\gk,\gd,B_{n+2})$.  Since
$C_{s}(\gk,\gd,B_{n+2}) \in \cF_{V_{\gk + \gw}}$, there is an
$f\map{V_{\gk + \gw}^{<\gw}}{V_{\gk + \gw}} $ such that
$s \cap \; \m{cl}_f \sub C_{s}(\gk,\gd,B_{n+2})$.
Since $C_{s}(\gk,\gd,B_{n+2}) \in a_{n+1}$, such an $f$ is in
$a_{n+1}$ so $a_{n+1} \cap V_{\gk + \gw} \in \m{cl}_f$.  Since
$s \in b$, we have (by property (2))
$b \cap V_{\gk + \gw} \in s$.
But $a_{n+1} \cap V_{\gk+\gw} = b \cap V_{\gk + \gw}$ so
$a_{n+1} \cap V_{\gk + \gw} \in C_{s}(\gk,\gd,B_{n+2})$.

	\end{enumerate}
	This completes the construction and the proof.

\newpage

\bibliographystyle{amsalpha}

\ifx\undefined\bysame
\newcommand{\bysame}{\leavevmode\hbox to3em{\hrulefill}\,}
\fi

\end{document}